# Two new topological indices based on graph adjacency matrix eigenvalues and eigenvectors

Juan Alberto Rodríguez-Velázquez[a,b] and Alexandru T. Balaban [b]

(a) Universitat Rovira i Virgili, Departament d'Enginyeria Informàtica i Matemàtiques Av. Països Catalans, 26 43007 Tarragona. Spain

(b) Texas A&M University at Galveston, Foundational Sciences, Galveston, TX 77553, United States

**Abstract**

The Estrada topological index EE, based on the eigenvalues of the adjacency matrix, is degenerate for cospectral graphs. By additionally considering the eigenvectors, two new topological indices are devised, which have reduced degeneracy for alkanes or cyclic graphs. Index $RV_a$ shows similarity to EE in ordering of alkanes with 8 to 10 carbon atoms, whereas index $RV_b$ is more similar to the average distance-based connectivity (Balaban index J). Inter-correlations between these four topological indices are discussed, indicating which factors have predominant influence.

**Introduction**

Topological indices (TIs) are digital counterparts of chemical structures and therefore represent these discrete molecular constitutional formulas by numerical functions. Some of their main uses are for quantitative structure-property or structure-activity relationships, QSPR and QSAR, respectively [1-3]. Chemical graphs (representing atoms in molecules by vertices and covalent bonds by edges) are connected non-directed graphs which are graph-theoretically planar. Loops and multiple bonds can exceptionally be present if thus is specified.

Constitutional isomers of alkanes correspond to trees, and for hydrogen-depleted hydrocarbons, their vertex degrees are at most 4. The number of graph vertices is called the order of the graph, but since this article will be discussing the ordering of graphs induced by TIs, "order" will be seldom used in this numerical sense. When non-isomorphic graphs correspond to identical TIs, these are said to be degenerate. Degeneracy of TIs is a serious drawback. For most chemical graphs, stereoisomerism is ignored, as will be done in this article, therefore only constitutional (structural) isomers will be taken into account.

The earliest TIs such as Wiener index [4], the Zagreb indices [5], or the Hosoya index [6], which are integers, have a high degeneracy. By contrast, TIs that are non-integer numbers, such as molecular connectivity (Randić [7]), higher-order molecular connectivity (Kier-Hall [8,9]), or information-theoretic indices (Bonchev [10], Trinajstić [11]) have lower degeneracy. Another feature, concerning the nature of the graph-theoretic invariants that are at the basis of TIs, further reduces the degeneracy. On replacing local vertex invariants (LOVIs) based on the adjacency matrix (vertex degrees, i.e. sums over rows or columns) by distasums (i.e. sums over rows or columns in the distance matrix), one can obtain the average distance-based connectivity index J (Balaban index [12]):

$$J = \frac{m}{m-n+2} \sum_{ij \in E} \frac{1}{\sqrt{d(i)d(j)}}, \qquad (1)$$

where $E$ is the edge set of the graph, $m = |E|$ is the number of edges, $n$ is the number of vertices, and $d(i)$ is the sum of the distances from vertex $i$ to the remaining vertices of the graph (distasum). Furthermore, combining a triplet of LOVIs and operations, one can explore which combination reflects better the problem at hand [13,14]. One can devise a classification of TIs with decreasing degeneracy as 1$^{st}$ generation (integer TIs based on integer LOCIs), 2$^{nd}$ generation (real-number TIs based on integer LOVIs) and 3$^{rd}$ generation (real-number TIs based in real-number LOVIs) [2].

In general, a topological index can be expressed as a function $f(x_1, …, x_n)$ where $x_1, …, x_n$ are LOVIs. In this article, we discuss three cases in which the LOVIs are centrality measures and the function $f$ is the arithmetic mean, the quadratic mean or the weighted arithmetic mean.

**Topological indices obtained from centralities based on eigenvalues and eigenvectors**

In network analysis, it is important how a vertex is embedded in a relational network. Vertices that are in favorable structural positions often have more opportunities and fewer constraints than others. In this sense, an interesting kind of local characterization of networks is made numerically by using measures known as "centrality" [15]. There are several centrality measures that have been introduced and studied for real world networks, in particular for social networks. They account for the different vertex characteristics that allow them to be ranked in order of importance in the network.

There are various methods for defining and measuring centrality. For instance,
- the degree centrality assigns more weight to vertices of high degree;
- the closeness centrality is based on distasums and was the basis of several topological indices [16,17];
- the betweenness centrality depends on the number of shortest paths between pairs of vertices, passing through various vertices;
- the eigenvector centrality, which will be discussed in more detail below;
- The subgraph centrality, which will be also discussed in more detail below.

The eigenvector centrality, introduced in 1972 by Bonacich [18], assigns relative scores to all vertices in the network based on the principle that connections to high-scoring vertices contribute more to the score of the vertex in question than equal connections to low-scoring vertices. In order to introduce the eigenvector centrality we need some notation. Let $G = (V, E)$ be a connected graph, consisting of vertices $V$ and edges $E$. Let $A$ be the adjacency matrix of the graph; $A_{ij} = 1$ if vertices $i$ and $j$ are connected by an edge and $A_{ij} = 0$ if they are not. The following equation describes the eigenvector centrality $x_i$ of a vertex $i$.

$$x_i = \frac{1}{\lambda} \sum_{j \in V} A_{ij} x_j = \frac{1}{\lambda} \sum_{j \in V: ij \in E} x_j, \qquad (2)$$

where $\lambda$ is the largest eigenvalue of $A$ and $X = (x_1, …, x_n)^t$ is the positive eigenvector of $A$ corresponding to the eigenvalue $\lambda$. The centrality of a vertex is proportional to the sum of the centralities of the vertices to which it is connected. Consequently, a vertex has high value of eigenvector centrality either if it is connected to many other vertices or if it is connected to others that themselves have high centrality.

Another well-known centrality measure is the subgraph centrality which was introduced in 2005 by Estrada and Rodríguez-Velázquez [19]. The subgraph centrality characterizes the vertices in a network according to the number of closed walks starting and ending at the node. Closed walks are appropriately weighted such that their influence on the centrality decreases as the order of the walk increases. A walk of length $r$ is a sequence of (not necessarily different) vertices $v_1, v_2, …, v_{r+1}$ such that for each $i = 1, …, r$ there is a link from $v_i$ to $v_j$. A closed walk is a walk in which $v_{r+1} = v_1$. We denote the number of closed walks of length $k$ starting and ending on vertex $v_i$ by $\varpi_k(v_i)$. The subgraph centrality is defined as

$$s_i = \sum_{k=0}^{\infty} \frac{\varpi_k(v_i)}{k!}. \qquad (3)$$

As shown in [19], the subgraph centrality can be computed from the eigenvalues and eigenvectors of the adjacency matrix. Let $U_1, \ldots, U_n$ be an orthonormal basis of the Euclidean space $\mathbb{R}^n$ composed by eigenvectors of $A$ associated to the eigenvalues $\lambda_1, \ldots, \lambda_n$, respectively. Let $u_j^i$ denote the $i$-th component of $U_j$. The subgraph centrality may be expressed as follows:

$$s_i = \sum_{j=1}^{n} (u_j^i)^2 e^{\lambda_j}. \qquad (4)$$

The arithmetic mean of the subgraph centralities was proposed in [19] as a global structural measure of the network:

$$EE = \frac{1}{n} \sum_{i=1}^{n} S_i. \qquad (5)$$

Since $U_1, \ldots, U_n$ is an orthonormal basis of $\mathbb{R}^n$, we have that $\sum_{i=1}^{n} (u_j^i)^2 = 1$, for every $j = 1, \ldots, n$. Hence, the $EE$ index can expressed as

$$EE = \frac{1}{n} \sum_{i=1}^{n} e^{\lambda_i}. \qquad (6)$$

In fact, the EE index is proportional to the Estrada index, which is defined by $\sum_{i=1}^{n} e^{\lambda_i}$. The Estrada index was introduced in 2002 by Ernesto Estrada [20] as a measure of the degree of folding of a protein.

The Estrada index has the advantage that it can be calculated easily. Despite the fact that it can be obtained from the subgraph centralities, it only uses the information contained in the eigenvalues (it does not use the information contained in the eigenvectors) and, as a result, it does not discriminate between cospectral graphs, and there are several chemical compounds whose corresponding graphs are cospectral. For instance, in the Tables 2 and 3 it will be shown that 5 pairs of nonanes and two pairs of decanes are degenerate, where two compounds share the same value of the topological index.

Note that the graphs associated to all these chemical compounds are trees and, as stated in [21], for $n$ sufficiently large, almost every tree of $n$ vertices has a cospectral mate. Hence, it would be desirable to provide some topological index from the eigenvalues and eigenvector which would be able to distinguish between cospectral graphs. In this article we introduce two topological indices that might satisfy this requirement.

Our first new topological index ($RV_a$) uses all the information contained in the subgraph centralities of the vertices; our formula uses both the eigenvalues and the eigenvectors. It is defined as the square mean of the subgraph centralities of the vertices of the graph.

$$RV_a = \left( \frac{1}{n} \sum_{i=1}^{n} S_i^2 \right)^{1/2}. \qquad (7)$$

By definition, each vertex contributes to $RV_a$ according to its participation in all the closed walks of all lengths, where the participation in short closed walks contributes more than the participation in the large ones. For the studied cases, this topological index has shown no degeneracy.

Our second new topological index ($RV_b$) takes advantage of both the subgraph centrality and the eigenvector centrality. It is defined as the weighted arithmetic mean of the subgraph centralities where the weights are the eigenvector centralities of the vertices of the graph, which are appropriately normalized to satisfy $\sum_{i=1}^{n} x_i = 1$:

$$RV_b = \sum_{i=1}^{n} x_i s_i. \qquad (8)$$

The contribution of the subgraph centrality of the vertices to $RV_b$ is pondered on the principle that connections to high-scoring vertices contribute more to the score of the vertex in question than equal connections to low-scoring vertices. In particular, if a graph is regular, then all the eigenvector centralities are equal and, in such a case, $RV_b$ becomes the Estrada index. For the studied cases, this topological index has shown no degeneracy.

**Constitutional isomers of alkanes up to decanes**

For obtaining adjacency matrices of all possible alkane isomers, the computer program MOLGEN was developed by Kerber and his coworkers [22,23] and it was used in the present article. Another program, also freely available, is Todeschini's DRAGON'[24].

One can expect practically all possible constitutional isomers of alkanes to be present in natural oil and gas deposits. Octane numbers constitute an experimentally accessible parameter, and they have been correlated with structures via other topological indices [25,26].

Table 1 presents all constitutional isomers of alkanes with 8 carbon atoms (octanes), together with the known Research Octane Numbers (RON). Tables 2 and 3 present all constitutional isomers of alkanes with 9 and 10 carbon atoms (nonanes and decanes). Four topological indices are also included: Estrada index (*EE*), Balaban index (*J*), and the two new TIs, $RV_a$ and $RV_b$. Degenerate *EE* values for nonane and decane isomers are indicated in boldface characters.

Table 1. All constitutional isomers of octanes and their TIs.

| Name | J-index Order | J-index Value | EE-index Order | EE-index Value | RVa-index Order | RVa-index Value | RVb-index Order | RVb-index Value | RON |
|---|---|---|---|---|---|---|---|---|---|
| Octane | 1 | 2.53006 | 1 | 2.09426 | 1 | 2.11445 | 1 | 2.18449 | -19 |
| 2-Methylheptane | 2 | 2.71584 | 2 | 2.11964 | 2 | 2.16539 | 2 | 2.25984 | 21.7 |
| 3-Methylheptane | 3 | 2.86207 | 3 | 2.12084 | 3 | 2.17058 | 3 | 2.29413 | 26.8 |
| 4-Methylheptane | 4 | 2.91961 | 4 | 2.12086 | 4 | 2.17087 | 4 | 2.30324 | 26.7 |
| 2.5-Dimethylhexane | 5 | 2.92782 | 6 | 2.14502 | 6 | 2.21518 | 5 | 2.31342 | 55.5 |
| 3-Ethylhexane | 6 | 3.07437 | 5 | 2.12207 | 5 | 2.17605 | 6 | 2.32775 | 33.5 |
| 2.4-Dimethylhexane | 7 | 3.09883 | 7 | 2.14624 | 7 | 2.22055 | 7 | 2.35149 | 65.5 |
| 2.2-Dimethylhexane | 8 | 3.11177 | 11 | 2.17288 | 11 | 2.27968 | 11 | 2.44868 | 72.5 |
| 2.3-Dimethylhexane | 9 | 3.17082 | 8 | 2.14747 | 8 | 2.22606 | 8 | 2.38365 | 71.3 |
| 3.4-Dimethylhexane | 10 | 3.29248 | 9 | 2.14867 | 9 | 2.23112 | 9 | 2.40195 | 76.3 |
| 3-Ethyl-2-Methylpentane | 11 | 3.35488 | 10 | 2.14870 | 10 | 2.23140 | 10 | 2.40626 | 87.3 |
| 3.3-Dimethylhexane | 12 | 3.37338 | 13 | 2.17534 | 13 | 2.29031 | 14 | 2.48754 | 75.5 |
| 2.2.4-Trimethylpentane | 13 | 3.38892 | 15 | 2.19831 | 15 | 2.32761 | 13 | 2.48688 | 100.0 |
| 2.3.4-Trimethylpentane | 14 | 3.46423 | 12 | 2.17411 | 12 | 2.28022 | 12 | 2.45028 | 102.7 |
| 3-Ethyl-3-Methylpentane | 15 | 3.58321 | 14 | 2.17777 | 14 | 2.30062 | 15 | 2.51639 | 80.8 |
| 2.2.3-Trimethylpentane | 16 | 3.62328 | 16 | 2.20200 | 16 | 2.34346 | 16 | 2.54387 | 109.6 |
| 2.3.3-Trimethylpentane | 17 | 3.70832 | 17 | 2.20323 | 17 | 2.34874 | 17 | 2.55890 | 106.1 |
| Tetramethylbutane | 18 | 4.02039 | 18 | 2.25664 | 18 | 2.45698 | 18 | 2.67015 | - |

Table 2. All constitutional isomers of nonanes and their TIs.

| Name | J-index Order | J-index Value | EE-index Order | EE-index Value | RVa-index Order | RVa-index Value | RVb-index Order | RVb-index Value |
|---|---|---|---|---|---|---|---|---|
| Nonane | 1 | 2.59508 | 1 | 2.11485 | 1 | 2.13343 | 1 | 2.20228 |
| 2-Methyloctane | 2 | 2.74669 | 2 | 2.13741 | 2 | 2.17837 | 2 | 2.26693 |
| 3-Methyloctane | 3 | 2.87662 | 3 | 2.13848 | 3 | 2.18296 | 3 | 2.30085 |
| 2.6-Dimethylheptane | 4 | 2.91466 | 7 | 2.15997 | 7 | 2.22242 | 4 | 2.30859 |
| 4-Methyloctane | 5 | 2.95482 | 4 | 2.13850 | 4 | 2.18323 | 5 | 2.31358 |
| 2.5-Dimethylheptane | 6 | 3.06082 | 8 | 2.16103 | 8 | 2.22692 | 7 | 2.34048 |
| 2.2-Dimethylheptane | 7 | 3.07299 | **16** | **2.18474** | 17 | 2.27967 | 17 | 2.44858 |
| 3-Ethyloptane | 8 | 3.09225 | 5 | 2.13957 | 5 | 2.18780 | 6 | 2.33670 |
| 2.4-Dimethylheptane | 9 | 3.15125 | 9 | 2.16108 | 9 | 2.22745 | 9 | 2.36372 |
| 2.3-Dimethylheptane | 10 | 3.15528 | **11** | **2.16215** | 12 | 2.23208 | 12 | 2.38852 |
| 4-Ethyloptane | 11 | 3.17534 | 6 | 2.13959 | 6 | 2.18807 | 8 | 2.34550 |
| 3.5-Dimethylheptane | 12 | 3.22305 | 10 | 2.16213 | 10 | 2.23168 | 10 | 2.37587 |
| 2.2.5-Trimethylhexane | 13 | 3.28071 | 24 | 2.20729 | 24 | 2.32181 | 18 | 2.46777 |
| 4-Ethyl-2-methylhexane | 14 | 3.30739 | **11** | **2.16215** | 11 | 2.23194 | 11 | 2.38325 |
| 3.4-Dimethylheptane | 15 | 3.32476 | 13 | 2.16324 | 13 | 2.23683 | 13 | 2.41245 |
| 3.3-Dimethylheptane | 16 | 3.33007 | **18** | **2.18692** | 20 | 2.28914 | 21 | 2.48971 |
| 2.3.5-Trimethylhexane | 17 | 3.37660 | **16** | **2.18474** | 16 | 2.27536 | 15 | 2.42474 |
| 3-Ethyl-2-methylhexane | 18 | 3.41009 | 14 | 2.16327 | 14 | 2.23709 | 14 | 2.41785 |
| 4.4-Dimethylheptane | 19 | 3.43105 | **20** | **2.18697** | 21 | 2.28965 | 23 | 2.49783 |
| 2.2.4-Trimethylhexane | 20 | 3.46726 | 25 | 2.20841 | 25 | 2.32665 | 22 | 2.49659 |
| 3-Ethyl-4-methylhexane | 21 | 3.49948 | 15 | 2.16434 | 15 | 2.24157 | 16 | 2.43267 |
| 2.3.4-Trimethylhexane | 22 | 3.57583 | **18** | **2.18692** | 18 | 2.28482 | 19 | 2.47057 |
| 2.4.4-Trimethylhexane | 23 | 3.57675 | 26 | 2.20953 | 26 | 2.33164 | 24 | 2.52111 |
| 2.2.3-Trimethylhexane | 24 | 3.58873 | 27 | 2.21065 | 27 | 2.33672 | 26 | 2.54713 |
| 3-Ethyl-3-methylhexane | 25 | 3.61739 | 22 | 2.18913 | 22 | 2.29883 | 25 | 2.52544 |
| 3-Ethyl-2.4-dimethylpentane | 26 | 3.67762 | **20** | **2.18697** | 19 | 2.28534 | 20 | 2.47932 |
| 2.3.3-Trimethylhexane | 27 | 3.70209 | **28** | **2.21176** | 29 | 2.34169 | 29 | 2.56541 |
| 2.2.4.4-Tetramethylpentane | 28 | 3.74642 | 34 | 2.25472 | 34 | 2.41818 | 31 | 2.58468 |
| 3-Ethyl-2.2-dimethylpentane | 29 | 3.79291 | **28** | **2.21176** | 28 | 2.34154 | 28 | 2.56265 |
| 3.3.4-Trimethylhexane | 30 | 3.80240 | 30 | 2.21283 | 30 | 2.34599 | 30 | 2.57436 |
| 3.3-Diethylpentane | 31 | 3.82468 | 23 | 2.19131 | 23 | 2.30824 | 27 | 2.55240 |
| 2.2.3.4-Tetramethylpentane | 32 | 3.87760 | 32 | 2.23435 | 32 | 2.38300 | 32 | 2.58871 |
| 3-Ethyl-2.3-dimethylpentane | 33 | 3.91921 | 31 | 2.21395 | 31 | 2.35094 | 33 | 2.58997 |
| 2.3.3.4-Tetramethylpentane | 34 | 4.01374 | 33 | 2.23658 | 33 | 2.39289 | 34 | 2.62350 |
| 2.2.3.3-Tetramethylpentane | 35 | 4.14473 | 35 | 2.26143 | 35 | 2.44735 | 35 | 2.69005 |

Table 3. All constitutional isomers of decanes and their TIs.

| Name | J-index | | EE-index | | RVa-index | | RVb-index | |
|---|---|---|---|---|---|---|---|---|
| | Order | Value | Order | Value | Order | Value | Order | Value |
| Decane | 1 | 2.6476048 | 1 | 2.1313243 | 1 | 2.1484895 | 1 | 2.2155292 |
| 2-Methylnonane | 2 | 2.7731889 | 2 | 2.1516268 | 2 | 2.1887066 | 2 | 2.2715343 |
| 3-Methylnonane | 3 | 2.8861628 | 3 | 2.1525872 | 3 | 2.1928143 | 3 | 2.3045160 |
| 2.7-Dimethyloctane | 4 | 2.9094720 | 9 | 2.1719293 | 9 | 2.2281979 | 4 | 2.3049634 |
| 4-Methylnonane | 5 | 2.9680150 | 4 | 2.1526093 | 4 | 2.1930539 | 5 | 2.3187566 |
| 5-Methylnonane | 6 | 2.9984191 | 5 | 2.1526096 | 5 | 2.1930631 | 6 | 2.3226997 |
| 2.6-Dimethyloctane | 7 | 3.0332963 | 10 | 2.1728897 | 10 | 2.2322333 | 7 | 2.3329711 |
| 2.2-Dimethyloctane | 8 | 3.0437584 | 28 | 2.1942203 | 28 | 2.2796586 | 29 | 2.4475733 |
| 3-Ethyloctane | 9 | 3.0869008 | 6 | 2.1535698 | 6 | 2.1971544 | 8 | 2.3408416 |
| 2.5-Dimethyloctane | 10 | 3.1244022 | 11 | 2.1729122 | 11 | 2.2324779 | 9 | 2.3521633 |
| 2.3-Dimethyloctane | 11 | 3.1296019 | **15** | **2.1738950** | 16 | 2.2368734 | 16 | 2.3901492 |
| 2.4-Dimethyloctane | 12 | 3.1600356 | 12 | 2.1729346 | 12 | 2.2327247 | 13 | 2.3687412 |
| 3.6-Dimethyloctane | 13 | 3.1681741 | 13 | 2.1738504 | 13 | 2.2362706 | 11 | 2.3603469 |
| 2.2.6-Trimethylheptane | 14 | 3.2054553 | 42 | 2.2145228 | 42 | 2.3176012 | 31 | 2.4571491 |
| 4-Ethyloctane | 15 | 3.2055351 | 7 | 2.1535922 | 7 | 2.1974027 | 10 | 2.3530739 |
| 2-Methyl-5-ethylheptane | 16 | 3.2555303 | 14 | 2.1738726 | 14 | 2.2365064 | 14 | 2.3702818 |
| 3.5-Dimethyloctane | 17 | 3.2685549 | **15** | **2.1738950** | 15 | 2.2367521 | 15 | 2.3852421 |
| 3.3-Dimethyloctane | 18 | 3.2769613 | 32 | 2.1961867 | 34 | 2.2881867 | 37 | 2.4894258 |
| 4-n-Propylheptane | 19 | 3.2950821 | 8 | 2.1536144 | 8 | 2.1976419 | 12 | 2.3613249 |
| 2.3.6-Trimethylheptane | 20 | 3.3014030 | 27 | 2.1941978 | 27 | 2.2755386 | 20 | 2.4108842 |
| 3.4-Dimethyloctane | 21 | 3.3088403 | 19 | 2.1748782 | 19 | 2.2411492 | 21 | 2.4160399 |
| 2.4.6-Trimethylheptane | 22 | 3.3374303 | 26 | 2.1932596 | 26 | 2.2716942 | 19 | 2.4033912 |
| 2.2.5-Trimethylheptane | 23 | 3.3555081 | 43 | 2.2154838 | 43 | 2.3214996 | 33 | 2.4723359 |
| 2-Methyl-3-ethylheptane | 24 | 3.3759294 | **20** | **2.1749004** | 20 | 2.2413839 | 23 | 2.4220362 |
| 2-Methyl-4-ethylheptane | 25 | 3.3907856 | 17 | 2.1739175 | 17 | 2.2369969 | 17 | 2.3946626 |
| 4.5-Dimethyloctane | 26 | 3.3977891 | **20** | **2.1749004** | 21 | 2.2413844 | 22 | 2.4217611 |
| 3-Methyl-5-ethylheptane | 27 | 3.4122569 | 18 | 2.1748557 | 18 | 2.2407822 | 18 | 2.4020957 |
| 4.4-Dimethyloctane | 28 | 3.4175118 | 34 | 2.1962316 | 36 | 2.2886672 | 41 | 2.4999928 |
| 2.3.5-Trimethylheptane | 29 | 3.4616741 | 29 | 2.1951810 | 29 | 2.2797417 | 25 | 2.4376438 |
| 2.5.5-Trimethylheptane | 30 | 3.4647270 | 45 | 2.2164898 | 45 | 2.3260086 | 42 | 2.5039116 |
| 2.2.4-Trimethylheptane | 31 | 3.4694656 | 44 | 2.2155513 | 44 | 2.3222135 | 40 | 2.4984887 |
| 4-Isopropylheptane | 32 | 3.4998573 | 22 | 2.1749228 | 22 | 2.2416281 | 24 | 2.4287843 |
| 2.4.5-Trimethylheptane | 33 | 3.5027177 | 30 | 2.1952035 | 30 | 2.2799844 | 28 | 2.4459290 |
| 2.2.3-Trimethylheptane | 34 | 3.5184258 | 50 | 2.2175407 | 51 | 2.3310803 | 52 | 2.5468503 |
| 4-Methyl-3-ethylheptane | 35 | 3.5299303 | 23 | 2.1758833 | 23 | 2.2456430 | 26 | 2.4410632 |
| 2.2.5.5-Tetramethylhexane | 36 | 3.5630018 | 69 | 2.2571176 | 69 | 2.4037182 | 50 | 2.5391176 |
| 3-Methyl-4-ethylheptane | 37 | 3.5636884 | 24 | 2.1758836 | 24 | 2.2456519 | 27 | 2.4422798 |
| 3-Methyl-3-ethylheptane | 38 | 3.5755049 | 39 | 2.1981758 | 39 | 2.2969351 | 45 | 2.5270563 |
| 2.3.4-Trimethylheptane | 39 | 3.5833271 | 33 | 2.1962092 | 32 | 2.2845405 | 34 | 2.4763120 |
| 2.5-Dimethyl-3-ethylhexane | 40 | 3.6033366 | 31 | 2.1952260 | 31 | 2.2802247 | 30 | 2.4529220 |
| 2.4.4-Trimethylheptane | 41 | 3.6256003 | 47 | 2.2165797 | 48 | 2.3269589 | 47 | 2.5291111 |
| 2.2-Dimethyl-4-ethylhexane | 42 | 3.6308443 | 46 | 2.2165123 | 46 | 2.3261050 | 43 | 2.5083439 |

| Compound | | | | | | | | |
|---|---|---|---|---|---|---|---|---|
| 2.3.3-Trimethylheptane | 43 | 3.6333919 | 53 | 2.2185467 | 53 | 2.3355720 | 56 | 2.5659188 |
| 3.3.5-Trimethylheptane | 44 | 3.6418625 | 48 | 2.2174958 | 50 | 2.3303699 | 46 | 2.5283882 |
| 2.2.4.5-Tetramethylhexane | 45 | 3.6842411 | 61 | 2.2368376 | 61 | 2.3636591 | 48 | 2.5292779 |
| 3.4.5-Trimethylheptane | 46 | 3.6854065 | 36 | 2.1971702 | 35 | 2.2884982 | 36 | 2.4880566 |
| 4-Methyl-4-ethylheptane | 47 | 3.6902938 | 40 | 2.1982205 | 40 | 2.2974050 | 49 | 2.5341882 |
| 3.4-Diethylhexane | 48 | 3.6982206 | 25 | 2.1768664 | 25 | 2.2499032 | 32 | 2.4590977 |
| 2-Methyl-3-isopropylhexane | 49 | 3.7280024 | 35 | 2.1962541 | 33 | 2.2850223 | 35 | 2.4872661 |
| 2.2.3.5-Tetramethylhexane | 50 | 3.7348237 | 62 | 2.2378663 | 62 | 2.3684513 | 55 | 2.5635860 |
| 2.3-Dimethyl-4-ethylhexane | 51 | 3.7561076 | 37 | 2.1971923 | 37 | 2.2887289 | 38 | 2.4917805 |
| 3.3.4-Trimethylheptane | 52 | 3.7783714 | 55 | 2.2195302 | 55 | 2.3396811 | 58 | 2.5769669 |
| 2.4-Dimethyl-3-ethylhexane | 53 | 3.7979077 | 38 | 2.1972148 | 38 | 2.2889704 | 39 | 2.4971640 |
| 2.4-Dimethyl-4-ethylhexane | 54 | 3.8025778 | 52 | 2.2185242 | 52 | 2.3351007 | 53 | 2.5535568 |
| 2.2-Dimethyl-3-ethylhexane | 55 | 3.8089258 | 54 | 2.2185691 | 54 | 2.3356687 | 57 | 2.5662121 |
| 2.3.4.5-Tetramethylhexane | 56 | 3.8139947 | 49 | 2.2175182 | 47 | 2.3269069 | 44 | 2.5198068 |
| 3.4.4-Trimethylheptane | 57 | 3.8231803 | 56 | 2.2195526 | 56 | 2.3399175 | 59 | 2.5803935 |
| 2.3.3.5-Tetramethylhexane | 58 | 3.8655589 | 63 | 2.2388951 | 63 | 2.3731169 | 60 | 2.5882000 |
| 3.3-Diethylhexane | 59 | 3.8747750 | 41 | 2.2001872 | 41 | 2.3058837 | 54 | 2.5602035 |
| 2.2.4.4-Tetramethylhexane | 60 | 3.8875947 | 70 | 2.2591757 | 70 | 2.4127893 | 65 | 2.6037615 |
| 2.2.3.4-Tetramethylhexane | 61 | 3.9418123 | 64 | 2.2398564 | 64 | 2.3769147 | 63 | 2.5967998 |
| 2.3-Dimethyl-3-ethylhexane | 62 | 3.9435702 | 58 | 2.2205586 | 58 | 2.3443942 | 62 | 2.5958708 |
| 3.3-Dimethyl-4-ethylhexane | 63 | 3.9711371 | 57 | 2.2205361 | 57 | 2.3440196 | 61 | 2.5906081 |
| 2.4-Dimethy-3-isopropylpentane | 64 | 3.9835002 | 51 | 2.2176081 | 49 | 2.3278580 | 51 | 2.5392798 |
| 3.4-Dimethyl-3-ethylhexane | 65 | 4.0204938 | 59 | 2.2215199 | 59 | 2.3482637 | 64 | 2.6034151 |
| 2.3.4.4-Tetramethylhexane | 66 | 4.0341182 | 66 | 2.2408627 | 66 | 2.3813325 | 67 | 2.6131435 |
| 2.2.4-Trimethyl-3-ethylpentane | 67 | 4.0729186 | 65 | 2.2399238 | 65 | 2.3776144 | 66 | 2.6068608 |
| 2.3.3.4-Tetramethylhexane | 68 | 4.0892890 | 67 | 2.2418914 | 67 | 2.3860915 | 69 | 2.6343423 |
| 2.2.3.3-Tetramethylhexane | 69 | 4.1017838 | 71 | 2.2632920 | 71 | 2.4315682 | 72 | 2.6918968 |
| 2-Methyl-3.3-diethylpentane | 70 | 4.1535139 | 60 | 2.2225481 | 60 | 2.3529529 | 68 | 2.6193066 |
| 2.3.4-Trimethyl-3-ethylpentane | 71 | 4.2289889 | 68 | 2.2429202 | 68 | 2.3907265 | 70 | 2.6499190 |
| 2.2.3.4.4-Pentamethylpentane | 72 | 4.2311322 | 74 | 2.2825666 | 74 | 2.4624380 | 71 | 2.6806227 |
| 3.3.4.4-Tetramethylhexane | 73 | 4.2817568 | 72 | 2.2652602 | 72 | 2.4396341 | 73 | 2.7083686 |
| 2.2.3-Trimethyl-3-ethylpentane | 74 | 4.3283429 | 73 | 2.2652827 | 73 | 2.4398626 | 74 | 2.7100771 |
| 2.2.3.3.4-Pentamethylpentane | 75 | 4.4038180 | 75 | 2.2856548 | 75 | 2.4763131 | 75 | 2.7335578 |

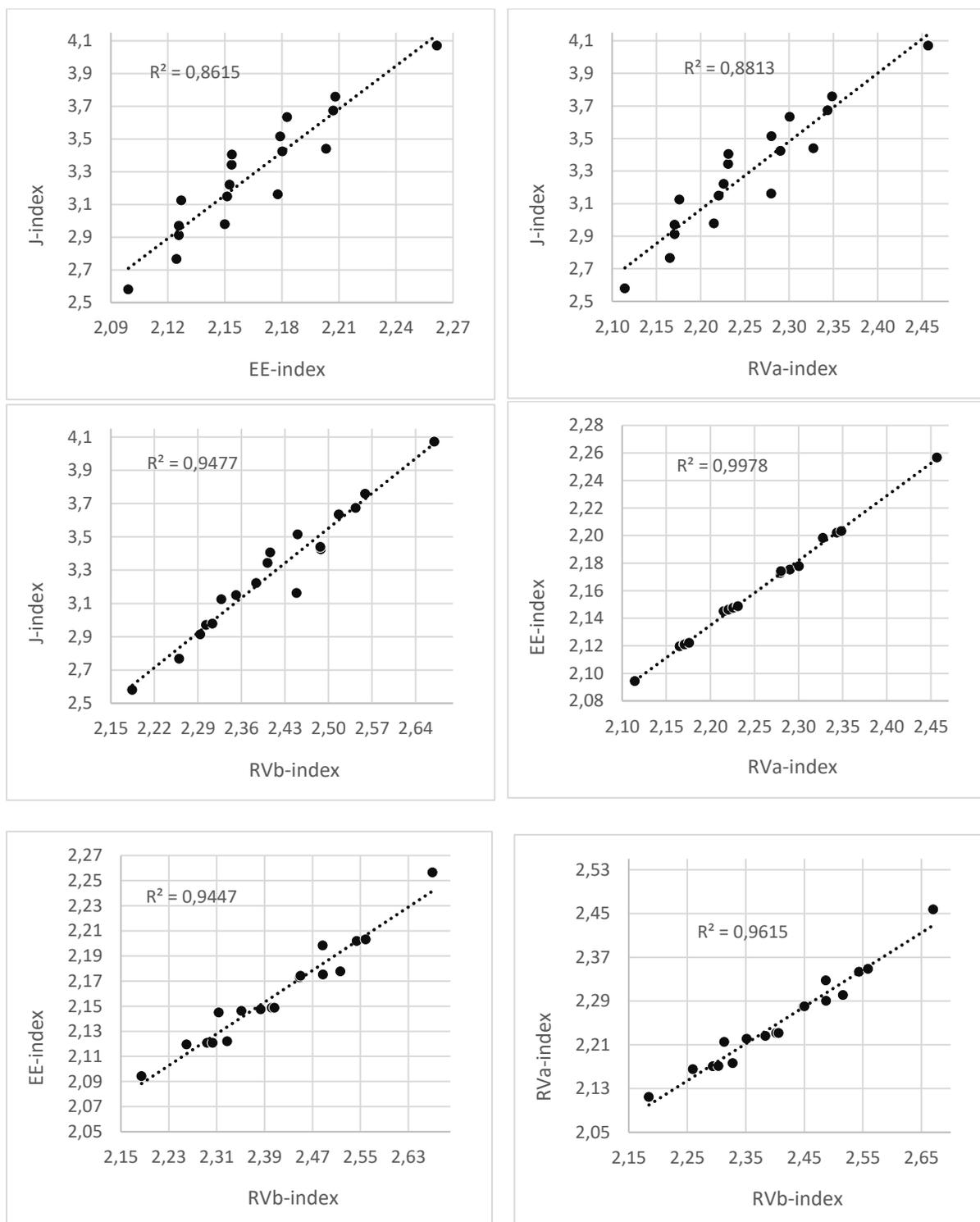

Figure 1. The six possible correlations between the four indices for constitutionals isomers of octanes.

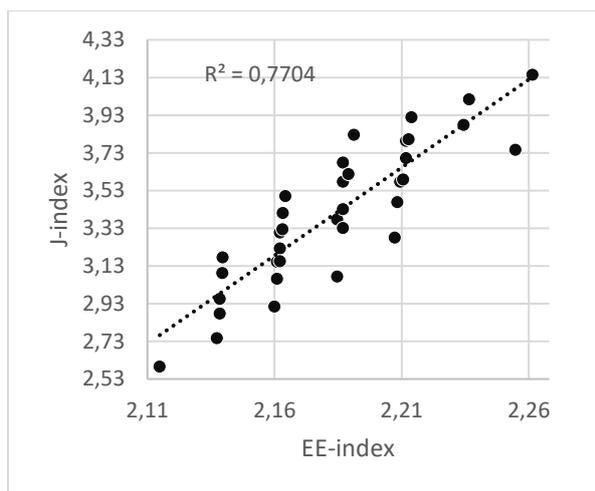
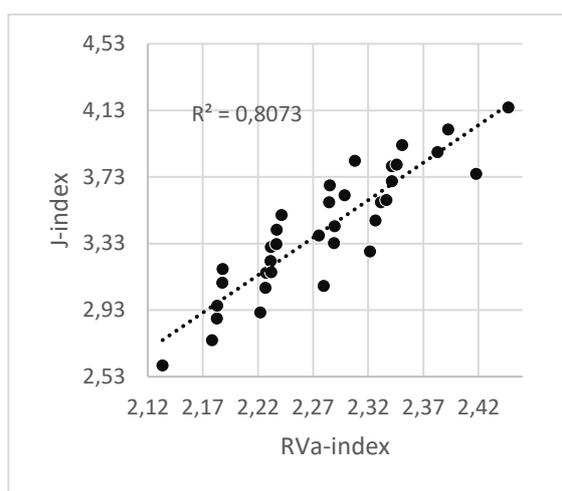
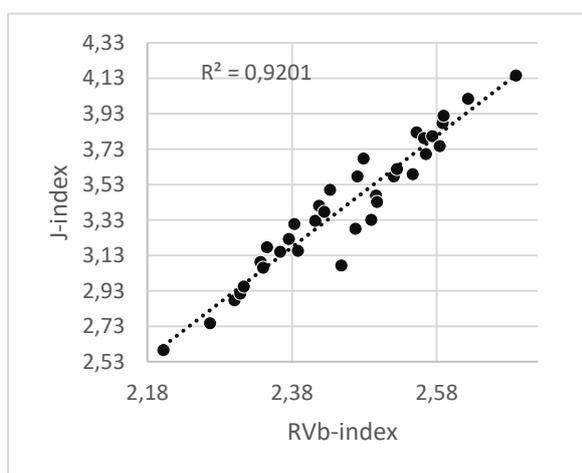
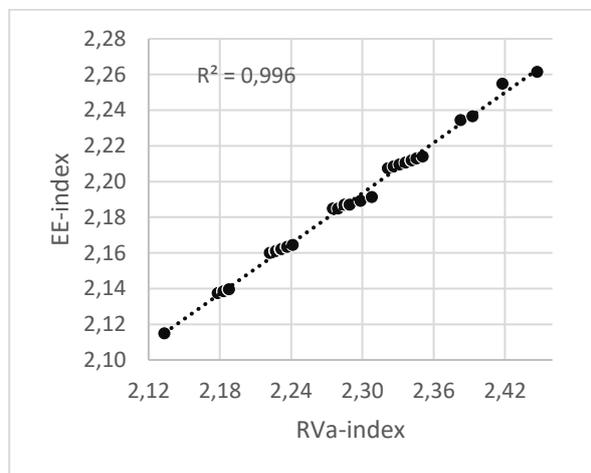
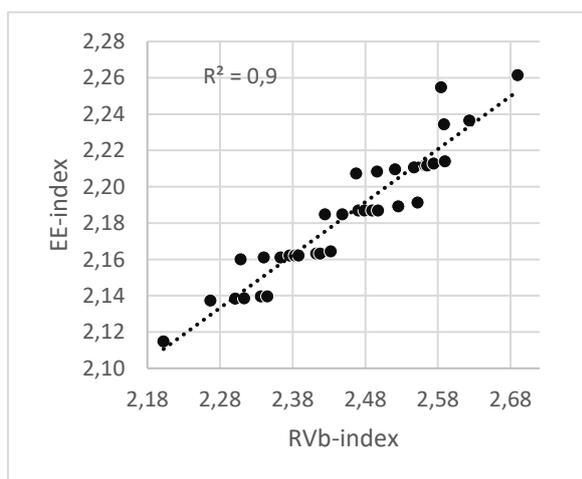
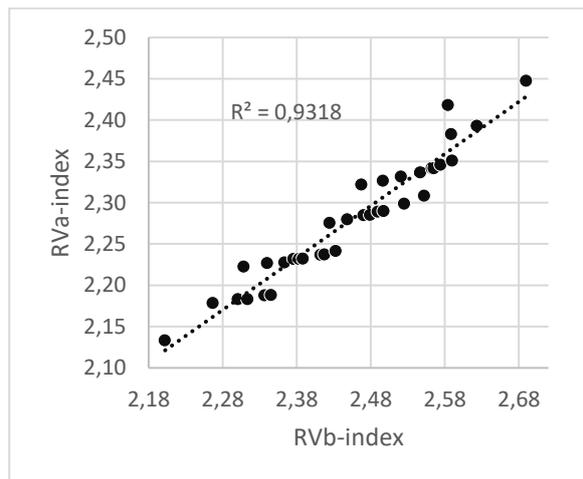

Figure 2. The six possible correlations between the four indices for constitutionals isomers of nonanes.

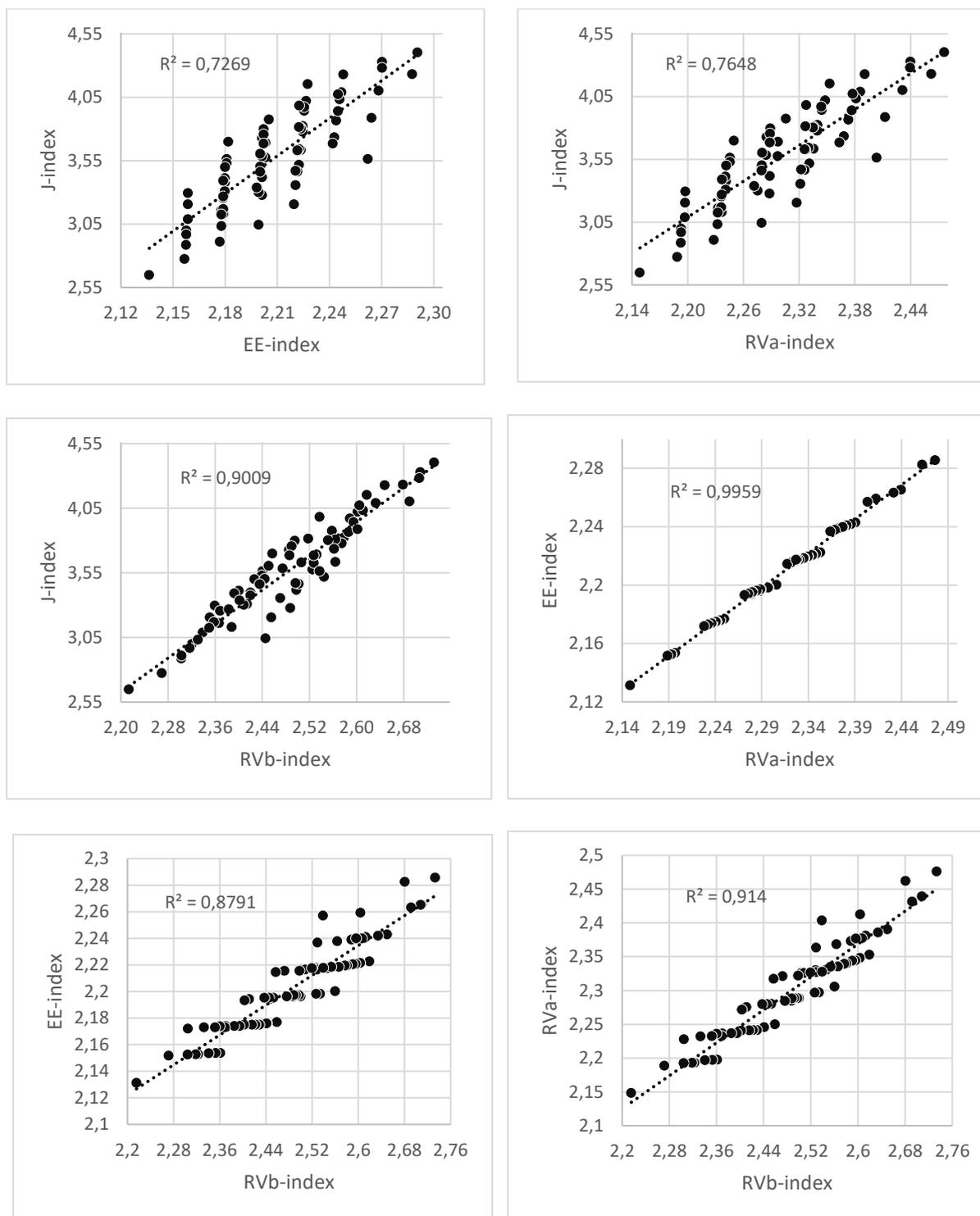

Figure 3. The six possible correlations between the four indices for constitutionals isomers of decanes.

Unlike most of the TIs of the 1st and 2nd generations, the four TIs examined in this article do not increase markedly with the size of the graph. Indeed, index $J$ for an infinitely long n-alkane has the asymptotic value of the number $\pi$. For the same number n of vertices representing carbon atoms, Figures 1, 2 and 3 show that in resulting clusters, the determining factors are number of vertices of highest degree, followed by the distance between them, and by the position of the high-degree vertices in the longest linear chain.

The tables for alkanes are organized as follows: Numerical values are limited to 7 significant digits. For each series of isomeric alkanes $C_nH_{2n+2}$, the corresponding alkane has a number starting with 1 in the column following the name, in the ordering dictated by index $J$ as seen in the next column. Then for indices $EE$, $RV_a$ and $RV_b$ the corresponding values are given, preceded by the ordering it would have if that index would dictated the order. For instance, 2,6-dimethylheptane has rank 4 according to $J$ and $RV_b$ indices, but rank 7 according to $EE$ and $RV_a$ indices.

The six pairwise inter-correlations between the four topological indices for alkanes with 8, 9, and 10 carbon atoms indicate an almost perfect agreement between $EE$ and $RV_a$ indices, and a fairly good agreement between $J$ and $RV_b$ The former pair has no overlapping values but a definite partition into non-overlapping clusters. Other inter-correlations also show clusters, but they have partial overlapping of values for one variable.

From previous discussions in the literature, the parameters that govern the alkane ordering are:
(a) Number of carbons atoms (order of graph)
(b) Longest linear chain or path (extremal disgrace)
(c) Number of vertices of highest degree
(d) Location of vertices of highest degree
(e) Distance between vertices of highest degree
(f) Number of and length of paths emerging from vertices of highest degree.

Interestingly, whereas the first and last alkane in each series agree for parameters (a) and (b), in between there are large variations among preferences for various parameters. From Figures 1 to 3 one can see that all four indices are well inter-correlated, and that $EE$ and $RV_a$ have the lowest overlap between clusters. One can predict the ordering for higher graph orders of trees.

Most substances have characteristic physical properties, and the easiest to determine are transition phase-transition temperatures, melting point and normal boiling point[1] (NBP). The former temperature (very close to the "triple point", the temperature at which all three phases coexist) is practically a constant for each substance, with little dependence on various external factors. The liquid state exists between the melting temperature and the critical temperature. The vapor pressure (and therefore, also the NBP) raises exponentially with increasing temperature, according to the Antoine equation or the Clausius-Clapeyon equation.

---

[1] One can easily identify the two types of scientists who call themselves "mathematical chemists". For those who have a theoretical background, the "boiling point" of a substance is a number taken from a Table in a textbook. On the other hand, those who handled substances in chemical laboratories and performed distillations under reduced pressure know and feel that" boiling point" is meaningless when the pressure conditions are not specified, therefore will always say or write "normal boiling point" indicating that it is the temperature at which the vapor pressure is equal to the normal pressure (760 Torr).

We found that none of the four indices examined in this article correlates with NBP, probably because neither graph-theoretical distances nor eigenvalues reflect intermolecular forces that govern phase transitions. Among experimentally measured properties of alkanes, we have chosen their behavior as energy source in internal combustion engines with spark ignition. With a scale from 0 for n-heptane and 100 for 2,2,4-trimethylpentane (iso-octane), one can compare the power and anti-knock properties of the alkane mixture (gasoline). In standard Otto-type engines the Research Octane Number (RON) is measured by maximum brake torquethe by comparison with binary mixtures of n-heptane and 2,2,4-trimethylpentane using the matching amount of iso-octane as the RON value. Structurally, the more branched alkanes have higher RON values than less branched ones. Improving RON values relies on catalytic cracking, partial aromatization and isomerization, and/or ethanol blending.

From Figure 4 one can see that satisfactory correlations exist between RON values and all four indices, and that *EE* leads to overlap-free clustering according to the presence of higher-degree vertices and then to the distance between them.

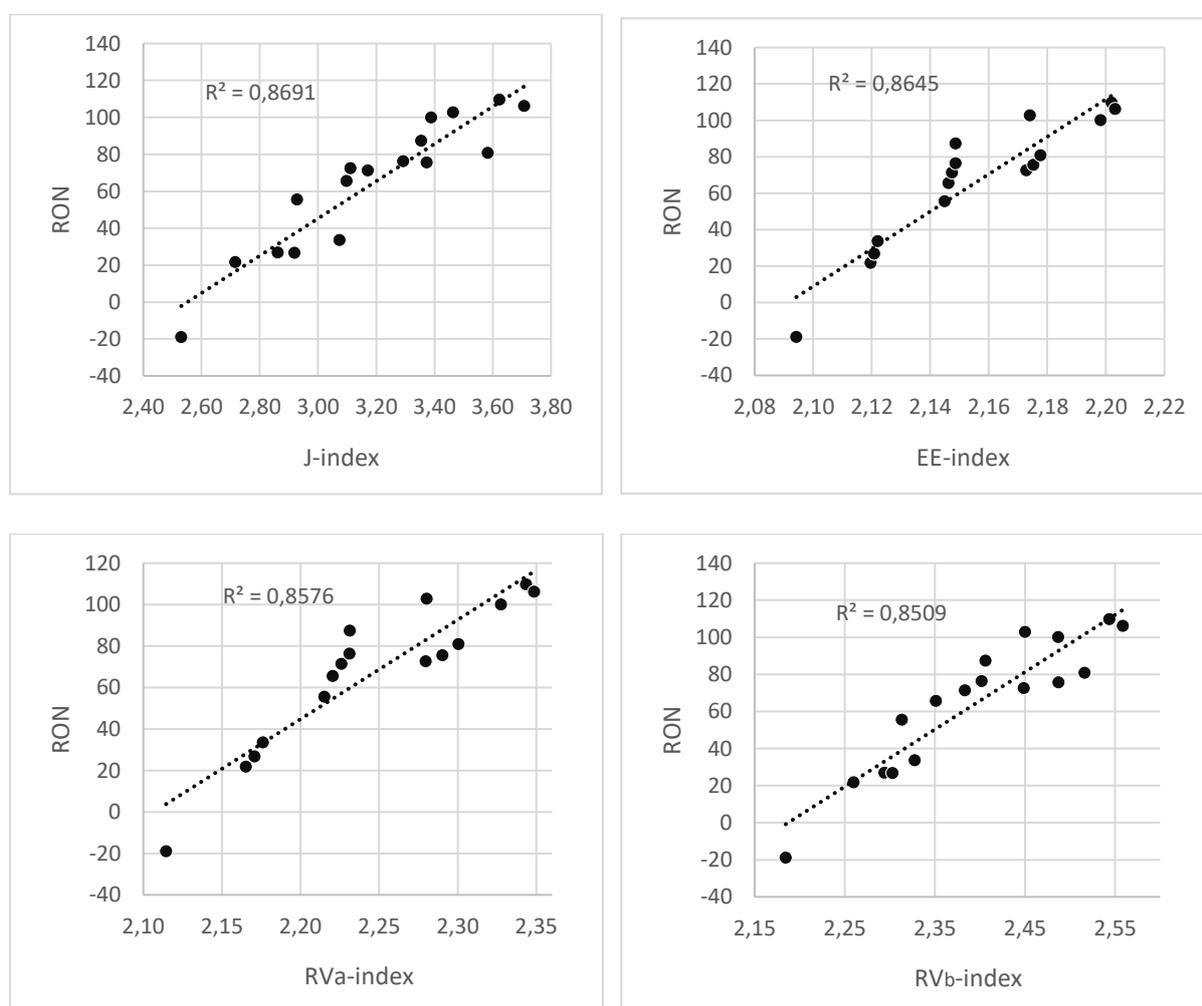

Figure 4. The correlation between the four indices for the octane and RON.

From the rich literature about how topological indices allow various orderings of alkanes, we mention only a few selected ones. Branching of threes can be influenced differently by various factors, as shown by Estrada, Rodríguez-Velázquez and Randić [27]. Randić [28], Trinajstić and coworkers [29], as well as Bonchev and coworkers [30], published tables comparing the ordering induced by various 1$^{st}$ and 2$^{nd}$ generation TIs. Balaban and coworkers [31,33] showed that the ordering of alkanes induced by index *J* differs from that induced by other TIs; since it parallels the ordering induced by the highly degenerate W index, it may be considered to be a "sharpened Wiener index" [34]. Bertz [35] used line graphs as the basis for his axioms and theorems, finding that the resulting ordering a surprising agreement with that determined by index *J*.

It should also be mentioned that before the *EE* index was introduced, several authors had used eigenvalues of chemical graphs: Lovasz and Pelican [36] had published six theorems for trees; Polanski and Giurmn [37] had commented on calculating the largest eigenvalue for molecular graphs; and Randić [28] had introduced two topological indices, $\lambda_1$ and $\lambda\lambda_1\lambda_1$.

**Cyclic chemical graphs**

All cyclic chemical graphs and vertex degrees from 1 to 4 with 5 and 6 vertices are presented in Figures 5 and 7, respectively. It should be mentioned that the 6-vertex graphs **32** and **35** were missing in ref. [32]. Owing to considerable steric strain, resulting from low bond angles and eclipsed conformations, most of the graphs with condensed cyclopropane and cyclobutane rings, representing the major part of Figures 5 and 7, are too unstable to exist in the real molecular world. Therefore these graphs have little chance to correspond to measurable experimental data. The discussion will involve the ordering of these chemical graphs.

On comparing the plots in Figure 8 among them, a marked difference is observed between the first three plots and the other ones. Such a difference is not observed in Figures 1, 2 and 3 related to acyclic graphs. This fact indicates that for acyclic graph the basic criteria for all four indices are similar, while for cyclic graphs there is a fundamental difference. We assume that apart from the topological distances, in the case of J index, the main selection criterion is the cyclomatic number, while for the other three indices the number of cycles and their size are also considered. In Figures 5 and 7, as well as Tables 4 and 5, the ordering of cyclic graphs will be based on index *EE*.

It is important to underline that "cycle" has different connotations in chemistry and in graph theory. Conventionally, the IUPAC-approved chemical nomenclature is based on the cyclomatic number μ. By definition, μ is the minimum number of edges one have to delete to convert the cycle structure into an acyclic one (a tree). The cyclomatic number is related to the number *n* of vertices and *m* of edges: $\mu = m - n + 1$. In chemistry, the number of cycles is determined by μ: although it is evident that decalin has two 6-membered rings and one 10-membered ring, it is named bicycle[4.4.0]decane.

Note that none of the four indices obeys exclusively the ordering determined by μ. In particular, index *EE* obeys this ordering for all graphs except graphs **44**, **47**, **61, 63** and **68**. Indices *EE*, $RV_a$ and $RV_b$ are clustered according to the actual number of closed walks: the smaller the walks, the greater its influence.

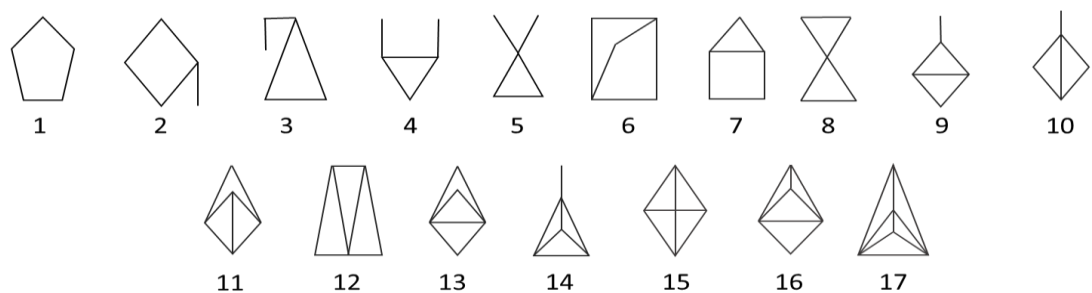

Figure 5. Cyclic chemical graphs with 5 vertices.

Table 4. Cyclic chemical graphs of 5 vertices and their TIs.

| G | μ | τ | Degrees | EE-index Order | Value | RVa-index Order | Value | RVb-index Order | Value | J-index Order | Value |
|---|---|---|---------|----------------|-------|-----------------|-------|-----------------|-------|---------------|-------|
| 1 | 1 | 0 | 2-2-2-2-2 | 1 | 2.29924 | 1 | 2.29924 | 1 | 2.29924 | 1 | 2.08333 |
| 2 | 1 | 0 | 3-2-2-2-1 | 2 | 2.40730 | 2 | 2.45286 | 2 | 2.50945 | 3 | 2.07967 |
| 3 | 1 | 1 | 3-2-2-2-1 | 3 | 2.60237 | 3 | 2.68055 | 3 | 2.83744 | 2 | 1.99894 |
| 4 | 1 | 1 | 3-3-2-2-1 | 4 | 2.66545 | 4 | 2.80036 | 4 | 2.96211 | 7 | 2.20344 |
| 5 | 1 | 1 | 4-2-2-1-1 | 5 | 2.70873 | 5 | 2.88230 | 6 | 3.02624 | 10 | 2.38220 |
| 6 | 2 | 0 | 3-3-2-2-2 | 6 | 2.93375 | 6 | 2.96019 | 5 | 2.97363 | 4 | 2.19089 |
| 7 | 2 | 1 | 3-3-2-2-2 | 7 | 3.07854 | 7 | 3.13230 | 7 | 3.17962 | 6 | 2.19393 |
| 8 | 2 | 2 | 4-2-2-2-2 | 8 | 3.32395 | 8 | 3.43256 | 8 | 3.49694 | 8 | 2.29966 |
| 9 | 2 | 2 | 3-3-3-2-1 | 9 | 3.43680 | 9 | 3.58056 | 9 | 3.74629 | 5 | 2.19235 |
| 10 | 2 | 2 | 4-3-2-2-1 | 10 | 3.50367 | 10 | 3.69958 | 10 | 3.84357 | 9 | 2.37197 |
| 11 | 3 | 2 | 3-3-3-3-2 | 11 | 4.04969 | 11 | 4.10869 | 11 | 4.14952 | 11 | 2.38901 |
| 12 | 3 | 3 | 4-3-3-2-2 | 12 | 4.34843 | 12 | 4.48283 | 12 | 4.54580 | 12 | 2.48607 |
| 13 | 3 | 3 | 4-4-2-2-2 | 13 | 4.51775 | 13 | 4.73073 | 13 | 4.80424 | 14 | 2.58080 |
| 14 | 3 | 4 | 4-3-3-3-1 | 14 | 4.87643 | 14 | 5.12857 | 14 | 5.35701 | 13 | 2.55465 |
| 15 | 4 | 4 | 4-3-3-3-3 | 15 | 5.57188 | 15 | 5.63395 | 15 | 5.64708 | 15 | 2.71108 |
| 16 | 4 | 5 | 4-4-3-3-2 | 16 | 6.02138 | 16 | 6.19395 | 16 | 6.26457 | 16 | 2.80428 |
| 17 | 5 | 7 | 4-4-3-3-3 | 17 | 8.04803 | 17 | 8.14364 | 17 | 8.16423 | 17 | 3.13746 |

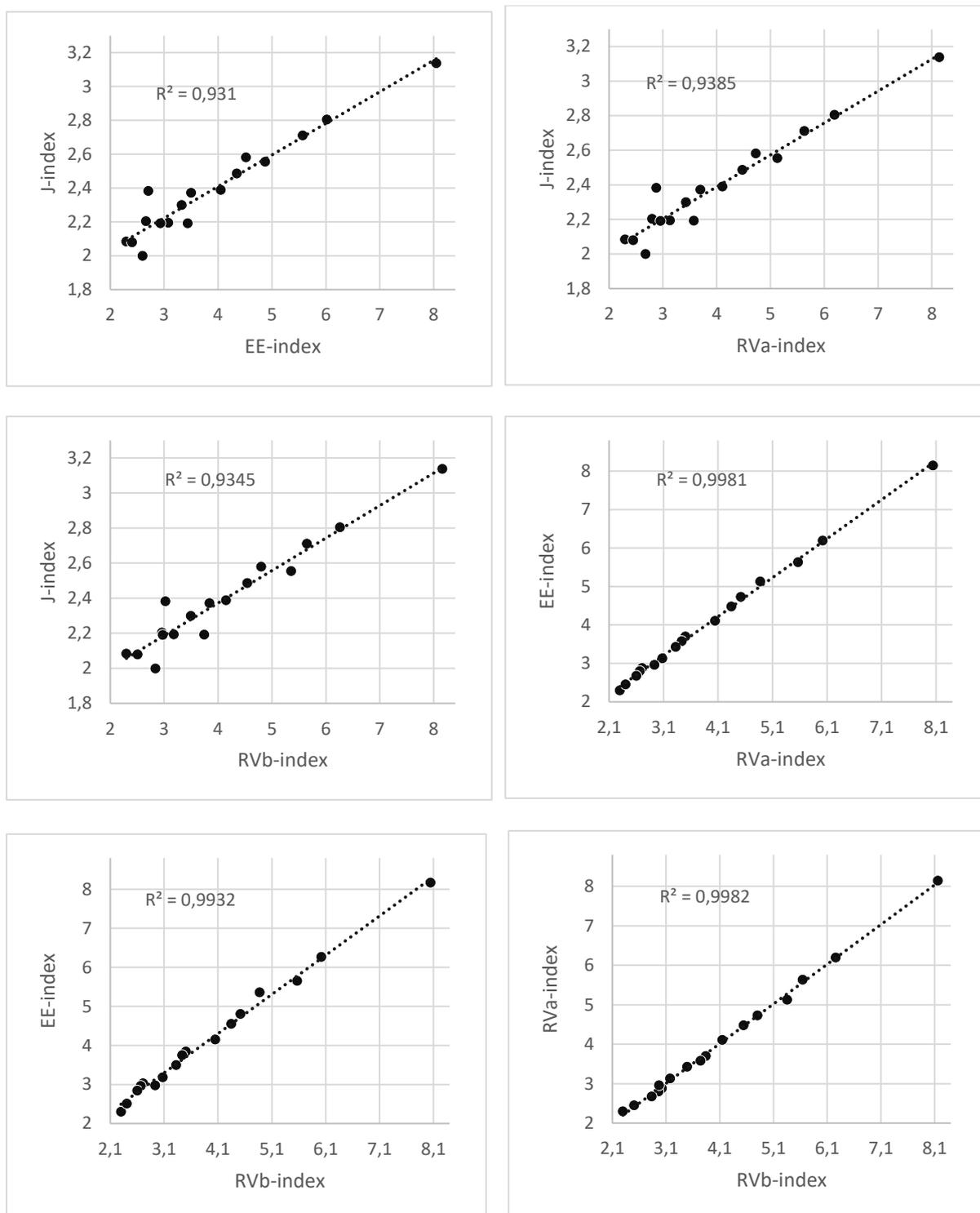

Figure 6. The six possible correlations between these four indices for the chemical graphs of 5 vertices.

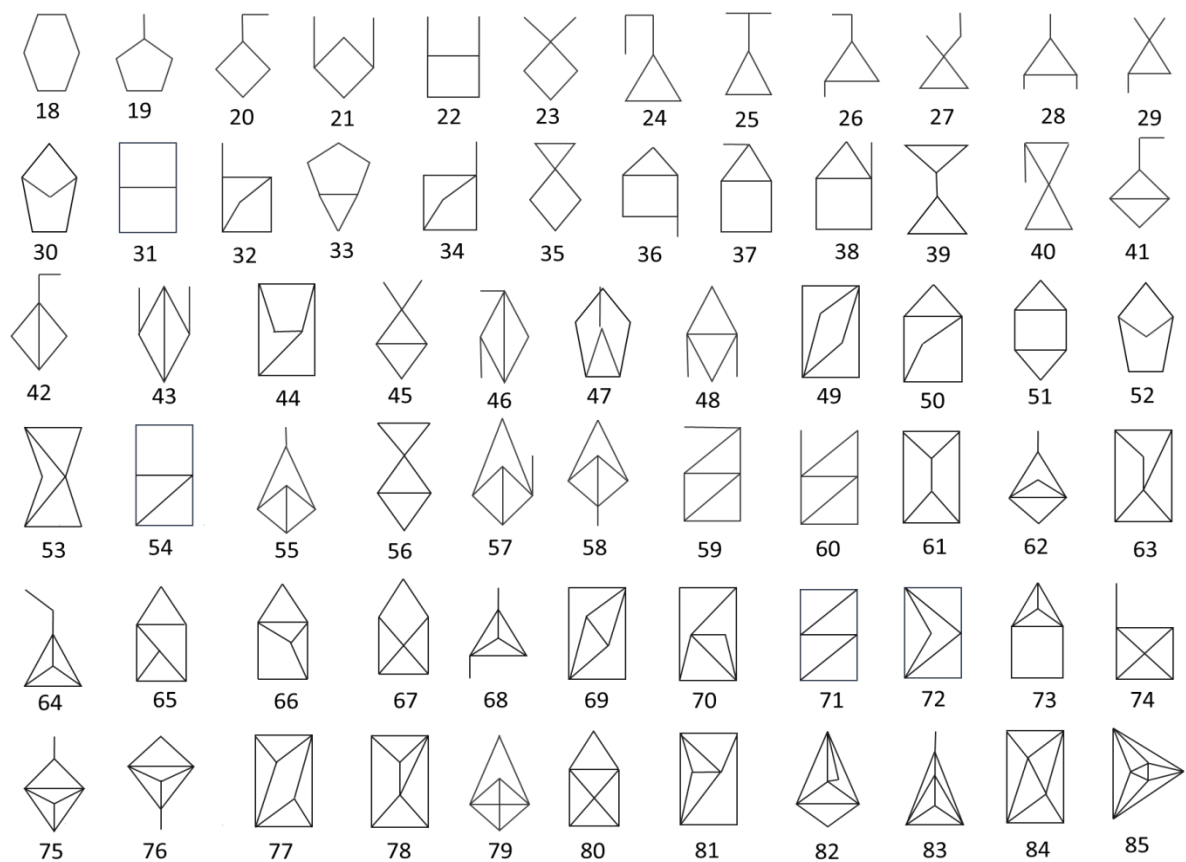

Figure 7. Cyclic chemical graphs of order 6.

Table 5. Cyclic chemical graphs of 6 vertices and their TIs.

| G | μ | τ | Degrees | EE-index Order | EE-index Value | RVa-index Order | RVa-index Value | RVb-index Order | RVb-index Value | J-index Order | J-index Value |
|---|---|---|---------|-------|-------|-------|-------|-------|-------|-------|-------|
| 18 | 1 | 0 | 2-2-2-2-2-2 | 1 | 2,28279 | 1 | 2,28279 | 1 | 2,28279 | 4 | 2,00000 |
| 19 | 1 | 0 | 3-2-2-2-2-1 | 2 | 2,33197 | 2 | 2,36690 | 2 | 2,41322 | 24 | 2,18411 |
| 20 | 1 | 0 | 3-2-2-2-2-1 | 3 | 2,38773 | 3 | 2,43202 | 3 | 2,53238 | 5 | 2,01427 |
| 21 | 1 | 0 | 3-3-2-2-1-1 | 4 | 2,42492 | 4 | 2,50077 | 4 | 2,60605 | 15 | 2,11237 |
| 22 | 1 | 0 | 3-3-2-2-1-1 | 5 | 2,42652 | 5 | 2,50704 | 5 | 2,61773 | 27 | 2,20160 |
| 23 | 1 | 0 | 4-2-2-2-1-1 | 6 | 2,46215 | 6 | 2,57847 | 6 | 2,68419 | 43 | 2,31144 |
| 24 | 1 | 1 | 3-2-2-2-2-1 | 7 | 2,54862 | 7 | 2,61848 | 8 | 2,82141 | 1 | 1,87629 |
| 25 | 1 | 1 | 3-3-2-2-1-1 | 8 | 2,58468 | 8 | 2,68305 | 10 | 2,85971 | 18 | 2,12947 |
| 26 | 1 | 1 | 3-3-2-2-1-1 | 9 | 2,60336 | 9 | 2,72974 | 11 | 2,95648 | 10 | 2,09391 |
| 27 | 1 | 1 | 4-2-2-2-1-1 | 10 | 2,64108 | 12 | 2,80605 | 13 | 3,02314 | 39 | 2,27633 |
| 28 | 1 | 1 | 3-3-3-1-1-1 | 11 | 2,65816 | 13 | 2,83751 | 15 | 3,06939 | 42 | 2,31136 |
| 29 | 1 | 1 | 4-3-2-1-1-1 | 12 | 2,69431 | 14 | 2,90619 | 17 | 3,12422 | 57 | 2,41329 |
| 30 | 2 | 0 | 3-3-2-2-2-2 | 13 | 2,71049 | 10 | 2,73622 | 7 | 2,75923 | 22 | 2,16250 |
| 31 | 2 | 0 | 3-3-2-2-2-2 | 14 | 2,75510 | 11 | 2,78823 | 9 | 2,82862 | 6 | 2,02774 |
| 32 | 2 | 0 | 3-3-3-2-2-1 | 15 | 2,87066 | 16 | 2,94549 | 14 | 3,04215 | 13 | 2,10846 |
| 33 | 2 | 1 | 3-3-2-2-2-2 | 16 | 2,87369 | 15 | 2,93174 | 12 | 3,01766 | 11 | 2,09486 |
| 34 | 2 | 0 | 4-3-2-2-2-1 | 17 | 2,90981 | 17 | 3,01702 | 16 | 3,10264 | 25 | 2,18660 |
| 35 | 2 | 1 | 4-2-2-2-2-2 | 18 | 2,96793 | 18 | 3,06429 | 19 | 3,15348 | 9 | 2,08964 |
| 36 | 2 | 1 | 3-3-3-2-2-1 | 19 | 2,98722 | 19 | 3,08006 | 20 | 3,20829 | 8 | 2,08348 |

| | | | | | | | | | | |
|---|---|---|---|---|---|---|---|---|---|---|
| 37 | 2 | 1 | 3-3-3-2-2-1 | 20 | 3,00322 | 21 | 3,11730 | 21 | 3,25996 | 16 | 2,11311 |
| 38 | 2 | 1 | 4-3-2-2-2-1 | 21 | 3,04446 | 22 | 3,19451 | 22 | 3,32489 | 32 | 2,21561 |
| 39 | 2 | 2 | 3-3-2-2-2-2 | 22 | 3,06779 | 20 | 3,09321 | 18 | 3,13569 | 2 | 1,91555 |
| 40 | 2 | 2 | 4-3-2-2-2-1 | 23 | 3,20784 | 24 | 3,37300 | 25 | 3,54063 | 23 | 2,16885 |
| 41 | 2 | 2 | 3-3-3-2-2-1 | 24 | 3,24630 | 25 | 3,40713 | 29 | 3,71379 | 3 | 1,92243 |
| 42 | 2 | 2 | 4-3-2-2-2-1 | 25 | 3,30426 | 28 | 3,51930 | 31 | 3,81147 | 12 | 2,09674 |
| 43 | 2 | 2 | 3-3-3-3-1-1 | 26 | 3,30592 | 27 | 3,51073 | 30 | 3,78694 | 7 | 2,05684 |
| 44 | 3 | 0 | 3-3-3-3-2-2 | 27 | 3,33162 | 23 | 3,36733 | 23 | 3,39946 | 20 | 2,15076 |
| 45 | 2 | 2 | 4-3-3-2-1-1 | 28 | 3,34446 | 29 | 3,57606 | 32 | 3,82336 | 31 | 2,21365 |
| 46 | 2 | 2 | 4-3-3-2-1-1 | 29 | 3,36395 | 33 | 3,62053 | 35 | 3,88217 | 33 | 2,23205 |
| 47 | 3 | 1 | 3-3-3-3-2-2 | 30 | 3,41566 | 26 | 3,47126 | 24 | 3,51917 | 29 | 2,21190 |
| 48 | 2 | 2 | 4-4-2-2-1-1 | 31 | 3,42039 | 35 | 3,72306 | 36 | 3,96981 | 46 | 2,33850 |
| 49 | 3 | 0 | 4-4-2-2-2-2 | 32 | 3,49632 | 31 | 3,60600 | 27 | 3,64775 | 40 | 2,30940 |
| 50 | 3 | 1 | 4-3-2-2-2-1 | 33 | 3,52030 | 32 | 3,61473 | 28 | 3,66747 | 36 | 2,26213 |
| 51 | 3 | 2 | 3-3-3-3-2-2 | 34 | 3,55616 | 30 | 3,58611 | 26 | 3,62031 | 21 | 2,15076 |
| 52 | 3 | 2 | 3-3-3-3-2-2- | 35 | 3,60141 | 34 | 3,70389 | 34 | 3,83449 | 30 | 2,21309 |
| 53 | 3 | 2 | 4-3-3-2-2-2 | 36 | 3,64171 | 36 | 3,75608 | 33 | 3,82659 | 37 | 2,26361 |
| 54 | 3 | 2 | 4-3-3-2-2-2 | 37 | 3,72074 | 37 | 3,87548 | 37 | 4,00570 | 28 | 2,20342 |
| 55 | 3 | 2 | 3-3-3-3-3-1 | 38 | 3,80332 | 38 | 3,94076 | 39 | 4,13722 | 17 | 2,11839 |
| 56 | 3 | 3 | 4-3-3-2-2-2 | 39 | 3,86059 | 39 | 3,97853 | 38 | 4,09360 | 19 | 2,14757 |
| 57 | 3 | 2 | 4-3-3-3-2-1 | 40 | 3,86567 | 40 | 4,04950 | 41 | 4,22145 | 34 | 2,23432 |
| 58 | 3 | 2 | 4-3-3-3-2-1 | 41 | 3,88241 | 41 | 4,08583 | 42 | 4,26077 | 35 | 2,25940 |
| 59 | 3 | 3 | 4-3-3-3-2-1 | 42 | 4,06971 | 43 | 4,30393 | 43 | 4,54374 | 26 | 2,19465 |
| 60 | 3 | 3 | 4-4-3-2-2-1 | 43 | 4,13279 | 44 | 4,41080 | 45 | 4,62346 | 41 | 2,31024 |
| 61 | 4 | 2 | 3-3-3-3-3-3 | 44 | 4,17908 | 42 | 4,17908 | 40 | 4,17908 | 44 | 2,31429 |
| 62 | 3 | 3 | 4-4-3-2-2-1 | 45 | 4,21678 | 46 | 4,53656 | 49 | 4,80111 | 38 | 2,26713 |
| 63 | 4 | 2 | 4-3-3-3-3-2 | 46 | 4,40006 | 45 | 4,50215 | 44 | 4,56136 | 49 | 2,36215 |
| 64 | 3 | 4 | 4-3-3-3-2-1 | 47 | 4,44911 | 50 | 4,78956 | 53 | 5,29320 | 14 | 2,10882 |
| 65 | 4 | 3 | 4-3-3-3-3-2 | 48 | 4,49826 | 47 | 4,57499 | 46 | 4,62867 | 50 | 2,36215 |
| 66 | 4 | 2 | 4-4-3-3-2-2 | 49 | 4,53316 | 48 | 4,67023 | 47 | 4,73992 | 54 | 2,40736 |
| 67 | 4 | 3 | 4-3-3-3-3-2 | 50 | 4,54190 | 49 | 4,67898 | 48 | 4,77110 | 51 | 2,36348 |
| 68 | 3 | 4 | 4-4-3-3-1-1 | 51 | 4,59935 | 53 | 5,02618 | 55 | 5,42737 | 52 | 2,36987 |
| 69 | 4 | 3 | 4-4-3-3-2-2 | 52 | 4,69034 | 51 | 4,87678 | 50 | 4,97489 | 55 | 2,41034 |
| 70 | 4 | 3 | 4-4-3-3-2-2 | 53 | 4,75921 | 52 | 4,98064 | 51 | 5,10382 | 47 | 2,35444 |
| 71 | 4 | 4 | 4-4-3-3-2-2 | 54 | 4,85595 | 54 | 5,04202 | 52 | 5,15816 | 48 | 2,35444 |
| 72 | 4 | 4 | 4-4-4-2-2-2 | 55 | 4,97186 | 55 | 5,23132 | 54 | 5,35595 | 59 | 2,45885 |
| 73 | 4 | 4 | 4-4-3-3-2-2 | 56 | 4,99092 | 56 | 5,27627 | 56 | 5,51432 | 56 | 2,41274 |
| 74 | 4 | 4 | 4-4-3-3-3-1 | 57 | 5,16652 | 57 | 5,44066 | 58 | 5,67859 | 53 | 2,38424 |
| 75 | 4 | 5 | 4-4-3-3-3-1 | 58 | 5,47348 | 60 | 5,83758 | 61 | 6,20599 | 45 | 2,32485 |
| 76 | 4 | 5 | 4-4-4-3-2-1 | 59 | 5,56433 | 61 | 5,97795 | 62 | 6,29889 | 58 | 2,45552 |
| 77 | 5 | 4 | 4-4-3-3-3-3 | 60 | 5,61537 | 58 | 5,65374 | 57 | 5,66970 | 60 | 2,53357 |
| 78 | 5 | 4 | 4-4-3-3-3-3 | 61 | 5,65728 | 59 | 5,74760 | 59 | 5,77892 | 61 | 2,53510 |
| 79 | 5 | 4 | 4-4-4-3-3-2 | 62 | 5,93489 | 62 | 6,11441 | 60 | 6,19802 | 63 | 2,57971 |
| 80 | 5 | 5 | 4-4-4-3-3-2 | 63 | 6,05872 | 63 | 6,24061 | 63 | 6,33070 | 64 | 2,58124 |
| 81 | 5 | 5 | 4-4-4-3-3-2 | 64 | 6,17775 | 64 | 6,44383 | 64 | 6,59387 | 65 | 2,58247 |
| 82 | 5 | 6 | 4-4-4-4-2-2 | 65 | 6,56292 | 65 | 6,89500 | 65 | 7,07465 | 62 | 2,57389 |
| 83 | 5 | 7 | 4-4-4-4-3-1 | 66 | 7,27213 | 67 | 7,73687 | 67 | 8,13602 | 66 | 2,60352 |

| 84 | 6 | 6 | 4-4-4-4-3-3 | 67 | 7,46309 | 66 | 7,55315 | 66 | 7,58039 | 67 | 2,76438 |
| 85 | 7 | 8 | 4-4-4-4-4-4 | 68 | 9,64480 | 68 | 9,64480 | 68 | 9,64480 | 68 | 3,00000 |

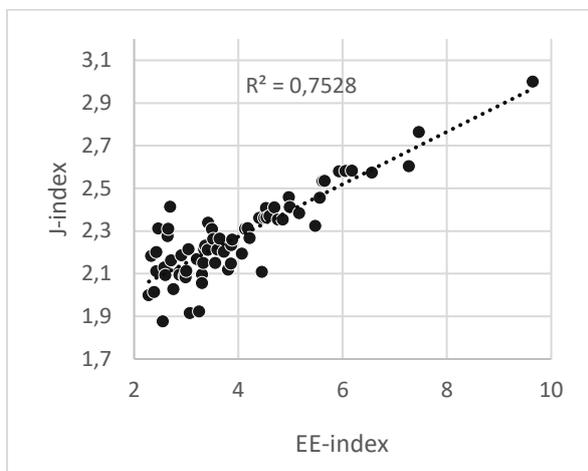
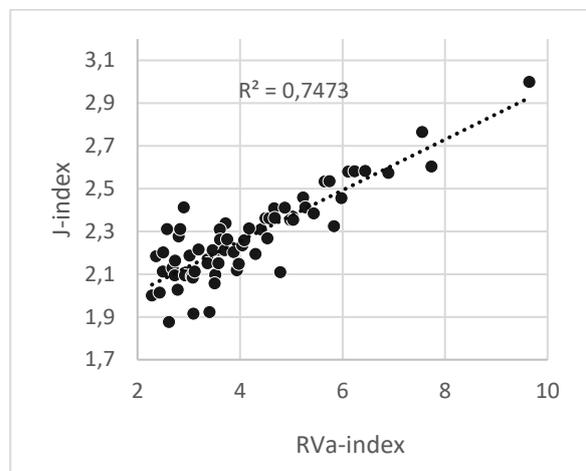
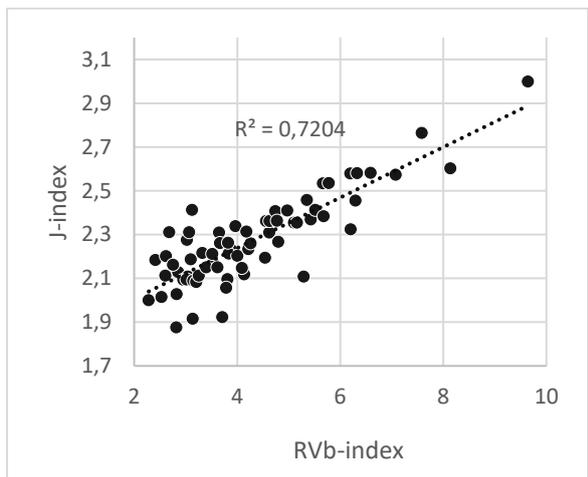
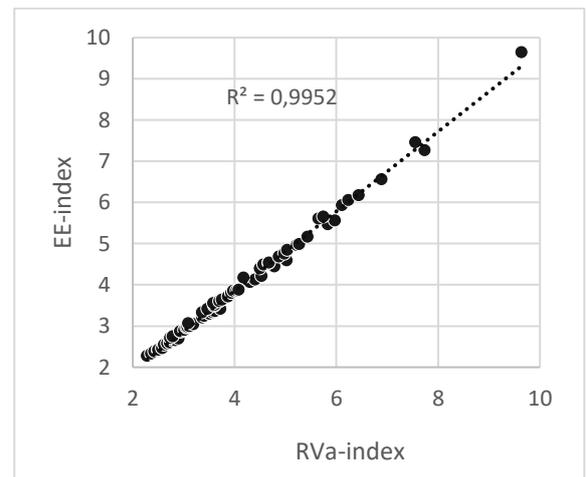
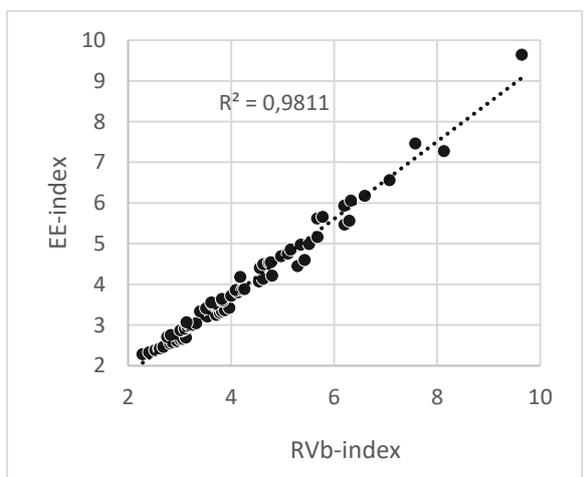
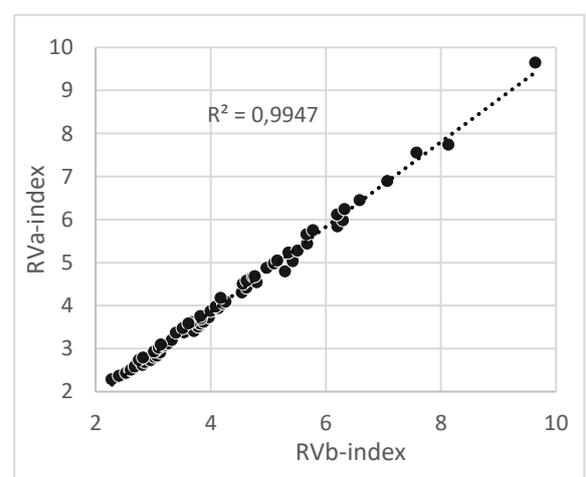

Figure 8. The six possible correlations between these four indices for the chemical graphs of 6 vertices.

Figures 6 and 8 show that for chemical graphs of 5 and 6 vertices the correlation between $EE$, $RV_a$ and $RV_b$ is almost perfect ($R^2 \geq 0{,}98$), while the correlation between $J$ and these three indices is lower ($R^2 \geq 0{,}72$). The rankings (in increasing order) determined by the different indices are shown in Tables 4 and 5. Note that none of these indices follow exclusively the order established by the cyclomatic number, while the number of small cycles (like triangles) and the degree sequence seems to be relevant in the order established by $RV_a$ and $RV_b$. Although $EE$ and $RV_a$ are very well correlated ($R^2 = 0{,}995$), there are some differences in the rankings imposed by them. For instance, graphs **29** and **30** are consecutive according to $EE$, while according to $RV_a$ the order is inverted, as graph **29** (which has one triangle) is in position 14 and graph **30** (which does not have any triangle) is in position 10. Similar behavior is observed when we compare graphs **83** (which has 7 triangles) and **84** (which has 6 triangles). Although these two indices are obtained from the subgraph centralities of the vertices, these differences in the order imposed by them might be a consequence of the information contained in Eq. 4 on the eigenvectors, which is no contained in $EE$ (see Eq. 6). The differences in the rankings imposed by $EE$ and $RV_b$ are more remarkable in the case of graphs **29** and **30**, which are in positions 12 and 13 according to $EE$, while they are in positions 17 and 7, respectively, according to $RV_b$. The small differences in the rankings established by $RV_a$ and $RV_b$ can be explained from the definition of both indices, as both are defined from the subgraph centralities but in the second one the values of the subgraph centralities are pondered by the values of the eigenvector centralities.

**Conclusions**

In order to avoid the degeneracy of topological indices based on eigenvalues of graphs (the Estrada index $EE$ becomes degenerative starting with trees having 9 vertices), we introduced two new topological indices, denoted by $RV_a$ and $RV_b$, which are based on eigenvalues and eigenvectors. Alkanes with 8, 9 and 10 vertices, denoting carbon atoms, and (poly) cyclic chemical graphs with 5 and 6 vertices were discussed in terms of inter-correlations and ordering according for four topological indices: $J$, $EE$, $RV_a$ and $RV_b$. A satisfactory correlation was obtained between Research Octane Numbers and the above four topological índices for alkanes of 8 carbon atoms.

This paper opens a challenge for mathematical chemists to search for the graph orders where degeneracy of the two new topological indices will first appear, as the number of graph vertices increases, for acyclic and for cyclic chemical graphs.

**Acknowledgement**

The first author of this paper thanks Prof. Douglas J. Klein, who invited him for a research stay at Texas A&M University at Galveston in 2018. The results included in this paper were obtained there. This research was supported in part by the Spanish government under the grants MTM2016-78227-C2-1-P and PRX17/00102.